\documentclass[review,12pt]{elsarticle}

\usepackage{hyperref}

\journal{Journal of Information Sciences}
\usepackage{savesym}
\savesymbol{AND}

\usepackage{amsfonts}
\usepackage{amsmath}
\usepackage{amssymb}
\usepackage{algorithm,algorithmic}
\usepackage{booktabs}
\usepackage{graphicx}
\usepackage{natbib}
\usepackage{color}
\usepackage{soul}
\usepackage[normalem]{ulem}
\usepackage{numcompress}

\makeatletter
\def\old@comma{,}
\catcode`\,=13
\def,{%
\ifmmode%
\old@comma\discretionary{}{}{}%
\else%
\old@comma%
\fi%
}
\makeatother

\usepackage{url}
\makeatletter
\g@addto@macro{\UrlBreaks}{\UrlOrds}
\makeatother

\newtheorem{theorem}{Theorem}

\newtheorem{remark}{Remark}

\newtheorem{definition}{Definition}
\newtheorem{assumption}{Assumption}

\newcommand{\ms}{{{\rm\bf q}}}

\renewcommand{\Pr}{{\mathbb{P}}}

\newcommand{\viol}{\texttt{viol}}

\newcommand{\decision}{{\mathbf{x}}}
\newcommand{\solvemip}{{\texttt{Solve\textsubscript{MIP}}}}
\newcommand{\verification}{{\texttt{Verification}}}
\newcommand{\feas}{{\texttt{feasible}}}
\newcommand{\seq}{{\texttt{seq}}}

\newcommand{\BB}{\mathcal{B}} 
 \newcommand{\FF}{\mathcal{F}}

\newcommand{\bQ}{{\mathbb{Q}}}

\newcommand{\bR}{{\mathbb{R}}}

\newcommand{\bZ}{{\mathbb{Z}}}

\renewcommand{\algorithmicrequire}{\textbf{Input:}}
\renewcommand{\algorithmicensure}{\textbf{Output:}}

\newcommand\oprocendsymbol{\hbox{$\square$}}
\newcommand\oprocend{\relax\ifmmode\else\unskip\hfill\fi\oprocendsymbol}

\begin{document}
\begin{frontmatter}
	\title{A Sequential Deep Learning Algorithm for Sampled Mixed-integer Optimisation Problems\tnoteref{mytitlenote}}
	
	\tnotetext[mytitlenote]{This work has been supported by the Australian Centre for
		Field Robotics and the Rio Tinto Centre for Mine Automation. A  preliminary  conference  version  of  this  paper  appeared  as  \cite{CHAMANBAZ20206749}.  This  paper  includes an entirely new learning-based algorithm, a more in-depth analysis, and completely new numerical computations. }
	
	
	\author[mohammadrezaaddress]{Mohammadreza  Chamanbaz\corref{mycorrespondingauthor}}
	\ead{m.chamanbaz@sydney.edu.au}
	\cortext[mycorrespondingauthor]{Corresponding author}
	\author[rolandaddress]{Roland Bouffanais} 
	\address[mohammadrezaaddress]{Rio Tinto Centre for Mine Automation, Australian Centre for Field Robotics,  the University of Sydney, Australia.}
	\address[rolandaddress]{Department of,
		Mechanical Engineering,
		Faculty of Engineering,
		University of Ottawa, Canada.}
	\ead{Roland.Bouffanais@uottawa.ca}

	\begin{abstract}
		 Mixed-integer optimisation problems can be computationally challenging. Here, we introduce and analyse two efficient algorithms with a specific sequential design that are aimed at dealing with sampled problems within this class. At each iteration step of both algorithms, we first test the feasibility of a given test solution for each and every constraint associated with the sampled optimisation at hand, while also identifying those constraints that are violated. Subsequently, an optimisation problem is constructed with a constraint set consisting of the current basis---namely, the smallest set of constraints that fully specifies the current test solution---as well as constraints related to a limited number of the identified violating samples. We show that both algorithms exhibit finite-time convergence towards the optimal solution. Algorithm \ref{alg:learning sequential algorithm} features a  neural network classifier that notably improves the computational performance compared to Algorithm \ref{alg:sequential algorithm}. We quantitatively establish these algorithms' efficacy through three numerical tests: robust optimal power flow, robust unit commitment, and robust random mixed-integer linear program. 
	\end{abstract}
	
	
\end{frontmatter}


\section{Introduction}
Sampled optimisation problems can comprise a great number of sampled constraints. They constitute an important category of problems that are widely encountered in the scenario approach~\citep{calafiore_uncertain_2004,calafiore2012mixed,CampiBook2018}, and also with learning paradigms at large, including statistical learning theory~\citep{vidyasagar_randomized_2001,vidyasagar_learning_2002,alamo_randomized_2009}. In these fields of study, semi-infinite optimisation problems---i.e., robust optimisation problems with an infinite number of constraints---are usually approximated and reformulated as a sampled optimisation problem associated with a finite set of random constraints. Effectively, the number of random samples is established so as to achieve a required level of probabilistic robustness when solving this approximate sampled optimisation problem. 
If a too stringent robustness criterion is sought, it is expected that the sample complexity---i.e., the size of the set of samples necessary to achieve this robustness criterion in a probabilistic senses---becomes high, thus yielding a computationally intensive sampled optimisation.

Dealing with the computational complexity of this particular class of optimisation problems has not attracted much attention despite its critical practical importance. In~\cite{chamanbaz_sequential_2013,chamanbaz_sequential_TAC_2016}, a sequential strategy has been considered in the frame of the scenario optimisation approach. This strategy consists in reducing at each iteration step $k$ the sample complexity, denoted as $N(k)$, as compared to the scenario bound. In a sequential fashion, this is followed by solving the sampled optimisation problem based on $N(k)$, and, ultimately by checking the robustness of the obtained solution by means of a validation test. The termination criterion is simply based on having the identified solution passing this test. Should the solution fail the validation test, a larger sample complexity $N(k+1)>N(k)$ is considered, and the previous steps are iterated until the algorithm ultimately satisfies the termination criterion. In \cite{chamanbaz_statistical_2014}, concepts from statistical learning theory are applied to solve robust linear and bilinear matrix inequality problems. Specifically, in this approach, the sample complexity is first computed and then a sequential randomised method is considered for the solution of the sampled optimisation problem. It is worth noting that the approaches reported in \cite{chamanbaz_sequential_2013,chamanbaz_statistical_2014,chamanbaz_sequential_TAC_2016} may end up demanding a large number of validation samples so as to satisfy the robustness condition of the candidate solution at each iteration. Furthermore, there are applications for which the set of samples is obtained from actual experiments, thereby making it a limited resource, and somehow restricting their availability. In \cite{calafiore2016repetitive}, a solution to the scenario problem is proposed based on concepts similar to those in~\cite{chamanbaz_sequential_TAC_2016}. Specifically, the algorithm introduced in~\cite{calafiore2016repetitive} does not rest on iterative increases of the cardinality of the set of design samples. Instead, it hinges on a probabilistic characterisation of the length of the iterative process required to arrive at a solution. Still using a sequential process, a `wait-and-judge scenario' optimisation \citep{Campi2018} has been developed without the need to test the validity of the candidate solution~\cite{GarattiIncremental2019}. In this algorithm, a sampled optimisation along with an estimation of the number of support constraints is iteratively carried out. At each iteration step, the robustness---in a probabilistic sense and based on set accuracy and confidence levels---controls the selection of the sample complexity. Practically, the method proposed in \cite{GarattiIncremental2019} primarily aims at minimising the number of scenarios without necessarily keeping the computational complexity low. Indeed, the repetitive evaluation of the number of support constraints may become computationally prohibitive when the size of the set of sample constraints grows.

Optimisation has long been at the core of Machine Learning (ML), and it continues to support the development of novel ML strategies. Interestingly, ML techniques have recently been considered to tackle demanding optimisation problems, such as in solving continuous and mixed-integer optimisation problems, see \cite{bengio2020machine} for an extensive review. Authors in \cite{lopez2016irace} used machine learning to tune the parameters of the optimisation algorithm automatically. Reinforcement learning \cite{sutton2018reinforcement} has been used in \cite{dai2017learning} to solve a diverse range of combinatorial problems defined over graphs where a neural network is trained to learn a heuristic algorithm which suggests the next node to visit contributing towards the optimal solution. In \cite{misra2018learning}, a statistical approach is developed to uncover the set of optimal active constraints---i.e., constraints that hold with equality in the optimal solution---for parametric optimisation problems. The algorithm needs to be fine-tuned for every problem and is limited to continuous optimisation problems---the paper claims that it can capture mixed-integer and non-convex problems but, no numerical evidence has been reported supporting that claim. The approach presented in \cite{bertsimas2019online,bertsimas2021voice} does not have the limitations of \cite{misra2018learning}. The authors of \cite{bertsimas2019online,bertsimas2021voice} consider parametric optimisation problems--- optimisation problems in which several parameters vary within a range each time the problem is formulated---and present a multiclass classifier to identify an \textit{optimal strategy} based on which one can recover the optimal solution without explicitly solving the problem. The optimal strategy is defined as the set of basis---also referred to as support constraints, see, e.g. \cite[Definition 2.1]{calafiore2010random}, which constitutes the minimal set of constraints underpinning the optimal solution---for continuous problems and the set of tight constraints together with the integer part of the optimal solution for mixed-integer problems. It is worth adding that none of the aforementioned approaches is designed to solve robust optimisation problems. 

In this paper---which is an extended version of \cite{CHAMANBAZ20206749}, we propose two algorithms for solving sampled optimisation problems. First, a sequential deterministic algorithm---Algorithm \ref{alg:sequential algorithm}---solves the sampled optimisation problem using concepts borrowed from the classical Las Vegas algorithm for linear and integer programming (LIP) \citep{Clarkson:1995:LVA:201019.201036}. Algorithm \ref{alg:sequential algorithm} as introduced in~\cite{CHAMANBAZ20206749} is intrinsically deterministic and does resort to any probabilistic validation at any stage. These features make it distinctive from other approaches, primarily seeking to maintain a low level of sampled constraints---as a way of keeping the computational burden for the optimisation task as reduced as possible---while maintaining probabilistic guarantees at the same level as those for the original problem. Here, we instead consider a computationally lean approach towards the solution of sampled optimisation problems for cases where the number of sampled constraints is large.  
The impetus for Algorithm \ref{alg:sequential algorithm} originates from a \textit{distributed} randomised constraints consensus approach introduced in \cite{chamanbaz_randomized_MILP_2017,chamanbaz2019randomized}, and that is effective at solving robust distributed mixed-integer problems. This specific approach is probabilistic---whereas both algorithms presented in the present paper are deterministic---with agents carrying out local computation and communication with the goal of achieving consensus on a candidate solution~\cite{chamanbaz_randomized_MILP_2017,chamanbaz2019randomized}. Note that a probabilistic validation step is included in \cite{chamanbaz_randomized_MILP_2017,chamanbaz2019randomized}, which has the same nature as the one used in~\cite{chamanbaz_sequential_TAC_2016}.

In Algorithm \ref{alg:learning sequential algorithm}, we borrow ideas from the classification technique used in \cite{bertsimas2021voice,bertsimas2019online} to further reduce the convergence time compared to Algorithm \ref{alg:sequential algorithm}. In particular, a neural network classifier is developed to predict the basis corresponding to  the optimisation problem formulated at an uncertain point. The classifier is used in reducing the computational complexity corresponding to the optimisation step of the algorithm. 

Given a sampled optimisation problem, the proposed two algorithms achieve a finite-time convergence towards the optimal solution. Both algorithms are sequential by design. Moreover, each iteration can be subdivided into the following steps: \emph{(i)} a validation step, and \emph{(ii)} an optimisation one. In that first step of this sampled optimisation, each and every constraint is verified for the candidate solution, and violating ones are singled out. In itself, stage \emph{(i)} is not computationally demanding since it is limited to validating the candidate solution without any optimisation \emph{per se}. In step \emph{(ii)}, the optimisation takes place with a set of constraints that comprises: \emph{(a)} the limited subset of constraints that have been identified in step \emph{(i)} as being violated, and \emph{(b)} the current \emph{basis}, which constitutes the smallest set of constraints fully characterising the test solution. These two steps---validation and optimisation ---are iteratively repeated until the test solution does not violate a single constraint, which is guaranteed to occur within a finite number of iterations. 

It is worth highlighting that the deterministic character of both algorithms enables their use for any optimisation problem exhibiting the Helly-type property; see \cite{amenta_helly-type_1994,amenta_hellys_2015}. 
As a matter of fact, this particular class of optimisation problems constitutes a sizeable share of non-convex problems. Lastly, it is worth stressing that a key novelty of this work---beyond the concepts reported in~\cite{CHAMANBAZ20206749} and the existing literature for this class of problems---is the use of machine learning (ML), and specifically a neural network multiclass classifier to solve robust optimisation problems. Compared to~\cite{CHAMANBAZ20206749}, the present paper contains an entirely new
ML-based algorithm, which seeks to identify participating constraints and thereby, exhibit a much reduced computational complexity.


\noindent
{\bf Notations\\}
In what follows, uppercase and italic letters, e.g., $F$, refer to constraints, whereas calligraphic uppercase letters, e.g., $\mathcal{F}$, refer to the set induced by the specific constraint $F$. Using these notations, if $A=B\cup C$ with $B$ and $C$ being collections of
constraints, then $\mathcal{A}=\mathcal{B}\cap \mathcal{C}$, that is, the set induced by the union of constraints $B$ and $C$ is the intersection of $\mathcal{B}$ and $\mathcal{C}$. Still with these notations, an optimisation problem
\begin{align*}
	\min\,\,& c^T\decision\\ 
	\text{subject to }& \decision\in \FF,
\end{align*}
is fully characterized by the pair $(F,c)$. Lastly, $J(F)$ is the smallest value of
$c^T\decision$ while $\decision\in\mathcal{F}$. Moreover, \texttt{randint(a,b)}  generates a random integer within the interval $[a,b]$ and $\lceil a\rceil$ returns the smallest integer number greater than $a$. 

\section{Preliminaries and Problem Statement}
Let us consider the following robust optimisation problem 
\begin{align}\nonumber
	\min\,\,& c^T\decision\\ \nonumber
	\text{subject to }& \decision\in \FF(q),\,\,\, \forall q\in\bQ \\ \label{eq:MICP}
	& \decision\in \bR^{d_R}\times\bZ^{d_Z},
\end{align}
where $\decision\in \bR^{d_R}\times\bZ^{d_Z}$ constitutes the vector of decision variables, the vector $q$ contains the entire set of uncertain parameters acting on the system, such that $\FF(q)=\{\decision\in \bR^d:f(\decision,q)\leq0\}$, with
$d=d_R+d_Z$ and, $f(\decision,q): \mathbb{R}^d\times\bQ\rightarrow\bR$,
is the constraint function of Problem \eqref{eq:MICP}. Note that when all decision variables $\decision$ are assumed to be continuous, then for any given value of $q$, this constraint $f(\decision,q)$ is a convex function. Without any loss of generality, a linear objective function may be assumed. Indeed, should a nonlinear convex objective function be considered, it could readily be transformed into the epigraph form by introducing an extra decision variable. 
Furthermore, if $d_Z=0$ then Problem \eqref{eq:MICP} reduces to a classical continuous convex optimisation task; if $d_R =0$, Problem \eqref{eq:MICP} constitutes an integer optimisation problem, while in the general case associated with $d_R\neq0,d_Z\neq 0$, we are in the presence of a mixed-integer optimisation. 

An efficient method to identify an approximate solution to \eqref{eq:MICP} with specific robustness guarantees---in the probabilistic sense---is recast it as a sampled optimisation problem by means of the scenario approach  \citep{calafiore_uncertain_2004,calafiore_scenario_2006,calafiore2010random,calafiore2012mixed,CampiBook2018}.
With this  approach, the semi-infinite optimisation problem \eqref{eq:MICP} is recast as a constrained optimisation with a finite number of constraints. These constraints are constructed based on some specifically samples, which are extracted randomly from uncertainty set $\bQ$. Formally, $N$ independent and identically distributed (i.i.d) samples are extracted from the set $\bQ$:
\[
\ms=\{q^{(1)},\ldots,q^{(N)}\}\in\bQ^N,
\]
and thus, one can formulate the following sampled optimisation problem
\begin{align}\nonumber
	\decision^*_N=\arg\min\,\,& c^T\decision\\ \nonumber
	\text{subject to }& \decision\in \bigcap_{i=1}^N\FF(q^{(i)})\\ \label{eq:sampled optimisation problem}
	& \decision\in \bR^{d_R}\times\bZ^{d_Z}.
\end{align}

It is important to note that only a few constraints are required to solve  Problem~\eqref{eq:sampled optimisation problem}. Indeed, the concept of \emph{basis}---the minimal set of constraints defining the solution---is key. This concept is closely connected to the Helly-type theorems originally proposed by E. Helly in \cite{helly1923mengen}, see \cite{amenta_helly-type_1994,amenta_hellys_2015} for more details. The central aim of this work is the identification of that basis for Problem \eqref{eq:sampled optimisation problem}, which offers a means to determine its solution $\decision^*_N$ in a computationally efficient fashion.   
\begin{definition}[Basis]
	Given a collection of constraints $F$, a subset of minimal cardinality $B\subseteq F$ is a basis of $F$ if the
	optimal cost of the problem defined by $(F,c)$ is identical to the one
	defined by $(B,c)$, and the optimal cost decreases if any constraint is removed.
\end{definition}
The \emph{combinatorial dimension} of any problem $(F,c)$ is defined by the size of its largest basis. Based on the following theorem adopted from  \cite[Corollary 1]{calafiore2012mixed},  and \cite[Theorem 3.11]{amenta_hellys_2015}, it can be explicitly expressed in terms of $d_R$ and $d_Z$ for the mixed-integer problem \eqref{eq:sampled optimisation problem}.
\begin{theorem}\label{thm:combinatorial dim}
	The combinatorial dimension of Problem \eqref{eq:sampled optimisation problem} is $d_{\text{comb}} = (d_R+1)2^{d_Z}-1$.
\end{theorem}
The following two assumptions are considered when seeking a solution to any subproblem of~\eqref{eq:sampled optimisation problem}. 
\begin{assumption}[Uniqueness]~\label{assum: uniqueness} For every $N\in\mathbb{N}$ and every multisample $\ms\in\mathbb{Q}^N$, there is a unique solution to Problem~\eqref{eq:sampled optimisation problem}.
\end{assumption}
Assumption~\ref{assum: uniqueness} is not \emph{per se} too restrictive. Indeed, to guarantee the uniqueness of the optimal point, one can resort to a number of strategies, such as considering a strictly-convex objective function, a lexicographic ordering, or a universal tie-breaking rule, see \cite[Observation 8.1]{amenta_helly-type_1994} for more details.
\begin{assumption}[Non-degeneracy]~\label{assum: degeneracy} The solution of any subproblem formed by every $N\in\mathbb{N}$ and every multisample $\ms\in\mathbb{Q}^N$ coincides with the solution of a problem which involves only basis constraints of the subproblem. 
\end{assumption}
Based on~\cite[Theorem 3 and Corollary 2]{calafiore2012mixed}, the robustness property of $\decision^*_N$ is fully characterised by the following theorem.
\begin{theorem}\label{thm:Scenario}
	Suppose Assumption \ref{assum: uniqueness} holds. Given probabilistic accuracy $\varepsilon\in(0,1)$ and confidence levels $\delta\in(0,1)$, let $N$ be the smallest integer satisfying 
	\begin{equation}\label{eq:scenario bound}
		\delta\geq \sum_{\ell=0}^{d_{\text{comb}}-1} {N \choose
			\ell} \varepsilon^\ell(1-\varepsilon)^{N-\ell}.
	\end{equation}
	Then the solution of \eqref{eq:sampled optimisation problem} $\decision^*_N$ satisfies
	\begin{align*}
		\Pr^{N}\bigg\{\ms\in\mathbb{Q}^{N}:
		\Pr\bigg\{q\in\mathbb{Q}:\decision^*_N\notin \FF(q)\bigg\}
		\leq\varepsilon \bigg\}\geq1-\delta,
	\end{align*}
	where $\Pr^N$ is the product probability measure on $\bQ^N$.
\end{theorem}
Note that~\eqref{eq:scenario bound} defines a binomial tail relating the number of samples $N$, with the accuracy level $\varepsilon$, confidence $\delta$, and the dimension of decision variables $d_R$ and $d_Z$. The sample complexity $N$ can be computed by numerically solving \eqref{eq:scenario bound} for any problems given desired accuracy and confidence levels.
Given that the number of scenario samples $N$ is inversely proportional to $\varepsilon$, and that there is a logarithmic relationship with $1/\delta$, so if $\varepsilon,\delta$ must somehow be made small, then $N$ in \eqref{eq:scenario bound} can become excessively large, thereby possibly yielding a challenging sampled optimisation. This possibility highlights the critical need for a computationally efficient algorithm towards the solution of the sampled optimisation problem \eqref{eq:sampled optimisation problem}. 

\section{Sequential Algorithm}
In this section, we detail the actual steps involved in our proposed sequential algorithm for the solution of sampled mixed-integer problems that can be formulated as in \eqref{eq:sampled optimisation problem}. As a first step, we introduce two fundamental primitives. The first, $[q^\texttt{Viol},\feas]=\verification(F(q^{(1)}),  \ldots,F(q^{(N)}),\decision,r) $, checks the feasibility of a test solution $\decision$  for all the sampled constraints involved in \eqref{eq:sampled optimisation problem}, i.e. $F(q^{(1)}),  \ldots,F(q^{(N)})$ and---if there exist---finds $r$ violating samples. In the event of a violation, the flag $\feas$ is set to $0$; otherwise, $\feas=1$. The second primitive is $[\decision,B]=\solvemip(F,c)$. This primitive deals with the actual solution of the optimisation problem constructed from $(F,c)$, and outputs the optimal point $\decision$ along with the associated basis $B$. The primitive needs first to solve the problem to find the optimal solution $\decision$, and then depending on the nature of the problem---being continuous, integer, or mixed-integer---it uses different methods to find the basis. See Remark \ref{rem:finding basis} for a detailed explanation of how to find the basis. These primitives form the backbone of the algorithm. First, a test solution $\decision$ is examined---using the $\verification$ primitive---to check if it is compatible with the full set of $N$ constraints of \eqref{eq:sampled optimisation problem} and---whenever possible---$r$ violating samples are identified; we note that there might be only $r'<r$ samples violating the candidate solution. Subsequently, the algorithm enters the optimisation phase---using the primitive $\solvemip$---whose constraint set consists of $(i)$ the constraints formed at $r$ violating samples and $(ii)$ the current set of basis. The algorithm then iterates over these two steps until no more violating samples can be found. The full process is formally detailed in  Algorithm~\ref{alg:sequential algorithm}.
\begin{algorithm}[t]
	\begin{algorithmic}[1]
		\caption{Sequential Algorithm}
		\label{alg:sequential algorithm}
		\STATE\algorithmicrequire{ $c,\,r,\,d_R,\,d_Z,\,F(q^{(i)}), i = 1, \ldots N,$}
		\STATE\algorithmicensure{ $\decision_\seq,B_\seq$}\\
		{\bf Initialisation:}
		\STATE Set $m=d_{\text{comb}} + 1$, $\feas=0$, $t=0$
		\STATE $[\decision(0),B(0)]=\solvemip(F(q^{(1)})\cup\ldots\cup F(q^{(m)}),c)$\\
		{\bf Evolution:}
		\WHILE{$\feas==0$}{
			\STATE $[q^\texttt{Viol},\feas]=\verification(F(q^{(1)}),  \ldots,F(q^{(N)}), \decision(t),r)$
			\STATE $[\decision(t+1),B(t+1)]=\solvemip(F(q^\texttt{Viol})\cup B(t),c)$
		}
		\STATE $t = t + 1$
		\ENDWHILE
		\STATE Set $\decision_\seq = \decision(t+1) $ and $B_\seq=B(t+1)$
		\RETURN $\decision_\seq,B_\seq$
	\end{algorithmic}
\end{algorithm}
Here are some important remarks related to Algorithm~\ref{alg:sequential algorithm}.
\begin{remark}[Complexity of an iteration step]
	A key attribute of Algorithm \ref{alg:sequential algorithm} is that the complexity of the optimisation task under consideration at each iteration step would not increase with the iteration count. Indeed, the \textit{maximum} possible number of constraints involved at the optimisation stage is $r+d_{\text{comb}}$. For instance, for a mixed-integer optimisation problem with $d_R = 5,\, d_Z = 3$ if we set $r = 10$, the number of constraints can be at most $57$; or in a continuous optimisation problem in which the dimension of the solution space is $d_R=5$, the number of constraints can be at most $15$. Therefore, at each iteration step, an optimisation of fixed---and small---complexity is solved. 
\end{remark}
\begin{remark}[Complexity of verification step]
	The verification step of Algorithm \ref{alg:sequential algorithm} is computationally inexpensive since it only requires checking the feasibility of a candidate solution for the $N$ sampled uncertainties present in Problem \eqref{eq:sampled optimisation problem}. 
\end{remark}
\begin{remark}[Identification of the basis]\label{rem:finding basis}
	In a continuous optimisation problem, the basis is fully characterised by the smallest set of active constraints, and it coincides with the set of active constraints if Assumption~\ref{assum: degeneracy} holds. Hence, in a continuous optimisation problem, it is straightforward to identify the basis. With a mixed-integer problem, however, identifying that basis can potentially become computationally demanding. Alternatively, a more tractable way of computing this basis---without necessarily seeking it to be of minimal cardinality---is to individually test each constraint to confirm whether or not they can belong to the basis.    
	In practice, one may consider dropping the $i$-th constraint for the optimisation sought. Should the objective value returned by this modified optimisation problem be smaller than the objective value of the original one, then the discarded constraint can be integrated into the basis. As previously highlighted, however, the number of constraints of the original problem---see line $7$ of Algorithm \ref{alg:sequential algorithm}---is at most $r+d_{\text{comb}}$. Thus, it is computationally inexpensive to identify the basis for problems with a limited number of constraints. We further note that a comparable greedy approach has been considered in the frame of algorithms aimed at solving scenario with discarded constraints \citep{campi2011sampling,calafiore2010random}, and more recently in the wait-and-judge scenario optimisation case \citep{Campi2018,GarattiIncremental2019}. In Section \ref{sec:learning}, we present a learning-based strategy to reduce the overall computational complexity further, while focusing on this particular step.

\end{remark}

\begin{remark}[Choice of the number of violating samples $r$]\label{rem:basis computation}
	The number of violating samples $r$ in Algorithm \ref{alg:sequential algorithm} is the byproduct of a trade-off between, on the one hand, the complexity of the optimisation task at each step and, on the other hand, the number of steps required for convergence. Specifically, a larger value for $r$ would yield a more complex optimisation at line $7$ of the algorithm but with a smaller iteration count. Depending on the computational capacity of the platform running Algorithm~\ref{alg:sequential algorithm}, one can tune $r$ to achieve the best performance in terms of total computational time.  
\end{remark}

The main features of Algorithm \ref{alg:sequential algorithm} are encapsulated by the following theorem. 
\begin{theorem}\label{thm:property of algorithm}
	Let Assumptions \ref{assum: uniqueness} and \ref{assum: degeneracy} hold. Then, the following statements hold. 
	\begin{enumerate}
		\item The objective value of the candidate solution $c^T\decision(t)=J(B(t))$ is monotonically increasing while Algorithm \ref{alg:sequential algorithm} is progressing. 
		\item 
		Algorithm \ref{alg:sequential algorithm} terminates in finite time.
		\item 
		The solution returned by Algorithm \ref{alg:sequential algorithm}, $\decision_\seq$ is identical to $\decision^*_N$.
	\end{enumerate}
\end{theorem} 
\vskip 1ex
{\bf Proof: } Note that when forming the basis at time $t+1$, i.e. $B(t+1)$, we use $B(t)$ from the previous time step $t$. Thus, $J(B(t+1))\geq J(B(t))$. Moreover, there has to be at least one violating constraint $F^\viol$ in all the iterations of Algorithm \ref{alg:sequential algorithm}, except of course at the last iteration. Indeed, if there had not been a violating constraint, the condition at line $5$ would have been satisfied and the algorithm would have terminated. This means that at line $7$ of Algorithm \ref{alg:sequential algorithm}, we solve an optimisation problem whose constraints set involves the current basis $B(t)$ and at least one violating constraint $F^\viol$. Therefore, owing to the presence of the violating constraint(s) $F^\viol$, and due to Assumption \ref{assum: uniqueness}, the cost has to increase, i.e. $J(B(t+1))> J(B(t))$. This completes the proof of statement (i) of Theorem \ref{thm:property of algorithm}. 

Since the number of constraints involved in \eqref{eq:sampled optimisation problem} is finite, then the number of candidate bases leading to a finite number of candidate costs $J(B(t))$ is also finite. Furthermore, as proved in the first statement, the cost $J(B(t))$ is strictly increasing with the iteration counter $t$. Since the sequence $\{J(B(t))\}_{t>0}$ is strictly increasing and has a finite number of elements, it will converge in a finite number of iterations leading to the finite-time termination of the algorithm. This completes the proof of statement (ii) of the theorem.  

We first note that since at any iteration $t$ of the algorithm, a subproblem of mixed-integer problem \eqref{eq:sampled optimisation problem} is being solved, $J(B(t))$ cannot be greater than $J(F)=c^T\decision^*_N$, where $F \doteq \bigcup_{i=1}^N F(q^{(i)})$; then, $J(B(t))\leq J(F),\,\forall t>0$ and as a result $J(B_\seq)\leq J(F)$. We now show that $J(B_\seq)$ cannot be smaller than $J(F)$. Assume by contradiction that $J(B_\seq)<J(F)$ or equivalently $J(B_\seq)<J(B_\seq\cup F)$ as $B_\seq\subseteq F$. By construction, $\decision_\seq$ must satisfy all the constraints in \eqref{eq:sampled optimisation problem} as it has passed the verification step of Algorithm \ref{alg:sequential algorithm}; then, $\decision_\seq\in\FF$ with $\FF\doteq \bigcap_{i=1}^N \FF(q^{(i)})$. Moreover, by definition  $\FF\subseteq \BB_\seq$, which implies that $\decision_\seq\in\FF\cap\BB_\seq$. Now, taking into account the fact that $\BB_\seq$ is the set generated by the basis $B_\seq$ associated with  $\decision_\seq$, we can state that $J(B_\seq)\geq J(F\cup B_\seq)$, which is in direct contradiction with our earlier assumption that $J(B_\seq)<J(B_\seq\cup F)$. Thus, $J(B_\seq)$ can neither be greater nor smaller than $J(F)=c^T\decision^*_N$. Therefore, $J(B_\seq)=c^T\decision^*_N$, which can be recast in the following equivalent form $c^T\decision_\seq=c^T\decision^*_N$. The latter, combined with Assumptions \ref{assum: uniqueness} and \ref{assum: degeneracy} concludes the proof of Theorem \ref{thm:property of algorithm}.

\section{Learning-based Sequential Algorithm}\label{sec:learning}

We first briefly discuss the approach first introduced in \cite{bertsimas2021voice}. Given a parametric optimisation problem with parameter $q$
\begin{align}\nonumber
	\min\,\,& c^T\decision\\ \nonumber
	\text{subject to }& \decision\in \FF(q),\,\,\,  \\ \label{eq:parametric_MILP}
	& \decision\in \bR^{d_R}\times\bZ^{d_Z},
\end{align}
an optimal strategy $S(q)$ is defined using which one can solve a reduced problem returning a solution identical to the optimal solution of \eqref{eq:parametric_MILP}. We remark that Problem \eqref{eq:parametric_MILP} is identical to \eqref{eq:MICP} when the set of uncertainty $\mathbb{Q}$ reduces to a singleton. The optimal strategy is defined as the set of basis for continuous problems and the group of active constraints jointly with the integer part of the decision variable at the optimal point for mixed-integer problems. For continuous problems, one can solve the problem only subject to the basis constraints, and for mixed-integer problems, the optimal value of the continuous part of the decision variable can be obtained by fixing the integer part to the value provided by the optimal strategy and solving a reduced problem whose constraint set involves only the set of active constraints. In \cite{bertsimas2021voice}, a multiclass classification problem is formulated to learn the mapping from parameter $q$ to the optimal strategy $S(q)$. A set of parameters $q_i\in\mathbb{Q}, i=1,\ldots,M$ is generated randomly, and for each parameter $q_i$, the optimal strategy $s_i$ is computed. The training data $(q_i,s_i), i=1,\ldots,M$, with $q_i$ as parameters and $s_i$ as the corresponding labels identifying the optimal strategy is used in training a multiclass classifier $\widehat{S}$. Given an unseen parameter $q_i$, the goal of the classifier $\widehat{S}$ is to identify  a strategy as close as possible to the optimal strategy. The classifier can assist in reducing the computational complexity of solving the parametric optimisation problem. The approach is very useful for parametric online optimisation problems where we repeatedly want to solve Problem \eqref{eq:parametric_MILP} for slightly different parameters $q$. In the online phase, where the goal is to solve \eqref{eq:parametric_MILP} for a particular parameter $q$, the parameter is first given to the trained classifier to estimate the optimal strategy, $\widehat{s} = \widehat{S}(q)$. Subsequently, the optimal strategy is used to find the optimal solution. Performing the mentioned two steps is much less computationally complex than directly solving \eqref{eq:parametric_MILP}.

The most computationally demanding part of Algorithm~\ref{alg:sequential algorithm} is to solve the mixed-integer problem and identify the corresponding basis at line 7. This step involves solving and finding the basis of an optimisation problem for which the set of constraints includes the current basis $B(t)$ along with the constraints formed by the violating samples $F(q^{\texttt{Viol}})$. The computational complexity of solving the problem and finding its basis has a direct relationship with the actual number of constraints involved in the problem. Hence, if we reduce the number of constraints without changing its solution or the set of basis, it is very likely that the complexity of the step presented at line 7 of Algorithm~\ref{alg:sequential algorithm} is reduced. This observation was the main motivation for using a learning-based strategy to reduce the computational complexity of Algorithm~\ref{alg:sequential algorithm}.
The computational burden can, in fact, be reduced using the approach presented in \cite{bertsimas2021voice}. To this end, we take a similar approach as \cite{bertsimas2021voice} to train a multiclass classifier $\widehat{S}$ which---having the violating sample $q^{\texttt{Viol}}$---can identify the basis constraint of an optimisation problem of the form \eqref{eq:parametric_MILP}. We use the classifier as an intermediate step to compute the basis of the problem $(F(q^\texttt{Viol})\cup B(t),c)$ (see line 7 of Algorithm \ref{alg:sequential algorithm}). The violating sample is first fed to the classifier $\widehat{S}$ to estimate  basis of the problem $B_{q^\texttt{Viol}} = \widehat{S}(F(q^\texttt{Viol}),c)$ and next, the estimated basis is used in $[\decision(t+1), B(t+1)]=\solvemip(B_{q^\texttt{Viol}}\cup B(t),c)$ to find basis of the problem $(F(q^\texttt{Viol})\cup B(t),c)$. It is worth noting that since the number of constraints in $B_{q^\texttt{Viol}}$ is much smaller than the ones  in $F(q^\texttt{Viol})$, the primitive $\solvemip(B_{q^\texttt{Viol}}\cup B(t),c)$ would be  computationally much cheaper than $\solvemip(F(q^\texttt{Viol})\cup B(t),c)$, see Tables \ref{tab:OPF results}, \ref{tab:unit commitment results}, and \ref{tab: simulation results MILP} which support this claim.

We remark that the optimal strategy used in \cite{bertsimas2021voice} for mixed-integer problems is the set of active constraints together with the integer part of decision variables at the optimal point. However, this strategy would not be useful in our case. For this reason, we modify the optimal strategy for the mixed-integer problem to be the basis, see Remark \ref{rem:basis computation} on how to compute the basis. To train the classifier, we first randomly generate several samples from the set of uncertainty $\mathbb{Q}$ and, for each sample, compute the basis. After generating all the training samples,  we use a one-hot encoding on the collection of the basis generated for all the training samples to define labels suitable for the multiclass classifier. There are a number of multiclass classifiers in the literature, however, we used a deep neural network to model multiclass classification \cite{bengio2007learning,lecun2015deep} as it shows superior performance compared to similar approaches such as support vector machine \cite{vapnik_statistical_1998} or random forest \cite{breiman2001random}. Section \ref{sec:simulation} provides more details on deep neural network architectures.

Any multiclass classifier has a certain probability of misclassification. In order to handle misclassification, we modify Algorithm~\ref{alg:sequential algorithm} so to minimise the effect or incorrect basis estimation. From Theorem \ref{thm:property of algorithm}, we know that the objective value of the candidate solution monotonically increases while Algorithm \ref{alg:sequential algorithm} progresses. A problematic misclassification is when a basis is incorrectly estimated, which results in the objective value of the candidate solution to stop increasing or even decrease. Such a problematic misclassification can  easily be recognised by comparing the objective value of the current candidate solution with the objective value of the candidate solution at the previous iteration. If the objective value is non-increasing, we use $[\decision(t+1),B(t+1)]=\solvemip(F(q^\texttt{Viol})\cup B(t),c)$ to find  the candidate basis and update the candidate solution. However, given that the classifier usually has very low probability of misclassification, the algorithm would rarely need to find the basis of the full problem. The full modified algorithm is detailed in Algorithm \ref{alg:learning sequential algorithm}. 
We limit the value of $r$---the number of violating samples returned by the verification primitive---to $1$ to simplify the classifier training procedure.

\begin{algorithm}[H]
	\begin{algorithmic}[1]
		\caption{Learning-based Sequential Algorithm}
		\label{alg:learning sequential algorithm}
		\STATE\algorithmicrequire{ $c,\,\widehat{S},\,d_R,\,d_Z,\,F(q^{(i)}), i = 1, \ldots N,$}
		\STATE\algorithmicensure{ $\decision_\seq,B_\seq$}\\
		{\bf Initialization:}
		\STATE Set $m=d_{\text{comb}} + 1$, $\feas=0$, $t=0$
		\STATE $[\decision(0),B(0)]=\solvemip(F(q^{(1)})\cup\ldots\cup F(q^{(m)}),c)$\\
		{\bf Evolution:}
		\WHILE{$\feas==0$}{
			\STATE $[q^\texttt{Viol},\feas]=\verification(F(q^{(1)}),  \ldots,F(q^{(N)}), \decision(t),1)$
			\STATE $B_{q^\texttt{Viol}} = \widehat{S}(q^\texttt{Viol})$
			\STATE $[\decision(t+1),B(t+1)]=\solvemip(B_{q^\texttt{Viol}}\cup B(t),c)$
			\IF{$J(B(t+1))\leq J(B(t))$}{
				\STATE 	$[\decision(t+1),B(t+1)]=\solvemip(F(q^\texttt{Viol})\cup B(t),c)$
			}
			\ENDIF
		}
		\STATE $t = t + 1$
		\ENDWHILE
		\STATE Set $\decision_\seq = \decision(t+1) $ and $B_\seq=B(t+1)$
		\RETURN $\decision_\seq,B_\seq$
	\end{algorithmic}
\end{algorithm}

\begin{remark}[Complexity of Algorithms \ref{alg:sequential algorithm}, and \ref{alg:learning sequential algorithm}]
    The complexity of Algorithms \ref{alg:sequential algorithm}, and \ref{alg:learning sequential algorithm} depends on two factors: (i) the specific type of problem being solved, and (ii) the algorithm used by the \solvemip~primitive. As mentioned earlier, formulation \eqref{eq:MICP} captures several broad classes of optimisation problems: convex continuous, integer and mixed-integer optimisation problems, including linear programming, mixed-integer linear programming, quadratic programming, mixed-integer quadratic programming, semi-definite programming, etc. Furthermore, several algorithms are available to be used by the \solvemip~primitive for any class of problems. This precludes a general and systematic complexity analysis of the two algorithms. In fact, the complexity analysis should only be considered for a pair problem-algorithm. For instance, linear programming problems can be solved using algorithms such as simplex, ellipsoid, and interior points. Hence, only for linear programming, several complexity analyses need to be presented, which is well beyond the scope of this paper. 

    A complexity analysis, however, can be presented in terms of the maximum number of iterations $t_{\max}$ required for Algorithms \ref{alg:sequential algorithm}, and \ref{alg:learning sequential algorithm} to solve the sampled optimisation problem \eqref{eq:sampled optimisation problem}. As stated in Theorem \ref{thm:combinatorial dim}, the combinatorial dimension of problem \ref{eq:sampled optimisation problem} is $d_{\text{comb}}$ which means the solution of \eqref{eq:sampled optimisation problem} is defined by at most $d_{\text{comb}}$ constraints. The first statement of Theorem \ref{thm:property of algorithm} proves that the objective value is monotonically increasing while Algorithm \ref{alg:sequential algorithm} is progressing. The same property is guaranteed for Algorithm \ref{alg:learning sequential algorithm} by the ``if" condition at line $9$. This excludes the possibility of repetition, that is, having two same basis sets in two different iterations of Algorithms \ref{alg:sequential algorithm}, and \ref{alg:learning sequential algorithm}. Subsequently, in the worst-case, Algorithms \ref{alg:sequential algorithm}, and \ref{alg:learning sequential algorithm} need to choose at most $d_{\text{comb}}$ out of all $N$ constraints. Therefore, the maximum number of iterations $t_{\max}$ in both Algorithms \ref{alg:sequential algorithm}, and \ref{alg:learning sequential algorithm} is bounded by
    \[
    t_{\max} \leq \max_{j = 1, \ldots, d_{\text{comb}}} {N \choose j }. 
    \]
\end{remark}


\section{Numerical Examples}\label{sec:simulation}
We have considered a comprehensive series of numerical simulations to thoroughly test the performance of Algorithms \ref{alg:sequential algorithm} and \ref{alg:learning sequential algorithm}. Specifically, we consider a wide range of different problems, including robust optimal power flow, robust unit commitment, and robust mixed-integer linear programming to assess the effectiveness of the proposed algorithms quantitatively. To this aim, the performance of these algorithms is compared in terms of the time required to complete the optimisation task. Furthermore, the performance of the presented algorithms is benchmarked against a direct solution of sampled optimisation problem obtained with widely used commercial solvers such as Gurobi \cite{gurobi} and Mosek \cite{andersen2000mosek}---we use Mosek for the robust optimal power flow problem which includes semi-definite constraints since Gurobi is unable to handle such constraints. 
All simulations are performed on a  Linux computing cluster in Rio Tinto Centre for Mine Automation. For all simulations, we allocated $12$ CPUs and $64$ GB of RAM. 

\subsection{Classifier Training}
Algorithm~\ref{alg:learning sequential algorithm} requires a classifier to estimate the basis at line 7. In this subsection, we discuss training the classifier and tuning its hyperparameters for all the numerical examples presented in the subsequent subsections. In order to train the classifier, we first need to generate several training samples $(q_i,s_i), i=1,\ldots, M$. The training sample is a tuple that includes uncertainty instance $q_i$ and its corresponding label $s_i$, which defines the optimal strategy. In order to generate training data, we start by extracting $M$ samples $q_i, i=1,\ldots,M$ from the uncertainty set $\mathbb{Q}$, and solve the parametric optimisation problem of the form \eqref{eq:parametric_MILP} formed at the extracted samples. Next, the basis is identified using the procedure mentioned in Remark \ref{rem:finding basis} and is encoded to form the optimal strategy. There are several ways to encode the optimal strategy; we use a one-hot encoding approach. This method creates a vector of all zeroes whose dimension is equal to the number of unique strategies found in the training data. If a data point belongs to the $i$th unique strategy, the $i$th component of the vector is set to $1$.

We need to train a  multiclass classifier to estimate the optimal strategy. We use a classical Neural Network (NN) approach to design the classifier. The neural network classifier has an input layer with a dimension equal to the number of uncertain parameters, several inner (hidden) layers, each with a depth that needs to be tuned, and an output layer whose dimension is equal to the number of unique strategies. The activation function is selected to be the rectified linear unit (ReLU) for all the layers except the last (output) layer, which has a softmax activation due to the multiclass nature of the classification problem.

We use the Keras library \cite{chollet2015keras} from  TensorFlow \cite{tensorflow2015-whitepaper} to implement the NN model. There are several hyperparameters---such as the number of hidden layers, depth of each hidden layer, batch size, number of epochs, and optimiser ---that should be tuned to design a classifier with the smallest misclassification error. Classically, $80\%$ of the data is used for training and the rest is used for testing the performance of the NN.  We use a grid search method and $K$-fold cross-validation from the scikit-learn library~\cite{scikit-learn}  to tune the hyperparameters. Moreover, we designed and tuned three classifiers for the problems discussed in the subsequent sections. The configuration and hyperparameters used for training the three NN classifiers are shown in Table~\ref{tab:classifier parameters}. We also report accuracy observed over the training and test sets. All  classifiers exhibit a very high accuracy on both training and test sets.

We use Matlab to model, generate the training samples, and solve the sampled optimisation problem for the optimal power flow and unit commitment problems while mixed-integer linear programming problem is modelled and solved in Python. For the optimal power flow and unit commitment problems the trained classifiers---which are trained using TensorFlow library---are exported to Matlab to predict the optimal strategy when using Algorithm~\ref{alg:learning sequential algorithm}.

\subsection{Robust Optimal Power Flow}
Optimal Power Flow (OPF) is an optimisation problem solved at regular intervals to define the operating point of controllable generators in power grids. Given the predicted demand, and network configuration, resources and limitations, OPF defines the active power of controllable generators and their magnitude of complex bus voltage so that the generation cost is minimised and the network constraints---such as line loading, min/max power rating of generators, and bus voltage---are respected. The increasing penetration of renewable energy resources introduces a large amount of uncertainty in the OPF. When uncontrollable resources fluctuate, the classical OPF solution can be very inefficient and may result in line overloads and potentially cascading outages. This calls for a robust strategy that generates policies that minimise the generation cost and, at the same time, ensure that the network constraints are not frequently violated.  

One of the successful approaches in designing a robust strategy is to use stochastic methods based on the scenario approach \cite{chamanbaz2019probabilistically,ETH-AC}, see \cite{chamanbaz2019Encyclopedia} for a full survey on available techniques. However, due to the complexity of the OPF problem and the fact that the number of decision variables is large, sampled optimisation Problem \eqref{eq:sampled optimisation problem} becomes very complex. For instance, for New England 39-bus systems case, choosing $\varepsilon=0.1, \delta = 1 \times10^{-10}$ the approach presented in \cite{ETH-AC} requires $11,831$ scenario samples and the sampled optimisation problem takes an impractical amount of time to get solved, see \cite[Section V]{chamanbaz2019probabilistically}.

\begin{table*}[t]
	\caption{Configuration and hyperparameters used for training neural network classifiers.}
	\begin{center}
		\scalebox{1}{
			\resizebox{1\textwidth}{!}{
				\begin{tabular}{c||c|c|c|c|c|c|c||c|c}
					\toprule
					Problem    & \# Training  & \# Unique&  \# Hidden  & \# Epochs & Batch & Optimiser &  Width of & Training    &   Test  \tabularnewline
					
					  &  Samples & Strategies & Layers & & Size & & Hidden Layers &     Accuracy    & Accuracy  \tabularnewline
					\midrule
					\midrule
					Optimal Power Flow   &  $28,000$  & $35$ &  $3$ & $200$ & $256$ & Adam &$512$ & $98.27\%$  & $98.25\%$     \tabularnewline \midrule
					Unit Commitment   &  $28,000$  & $138$ &  $3$ & $300$ & $128$ & Adam & $512$ & $96.58\%$  & $96.41\%$     \tabularnewline \midrule
					Mixed-integer Linear Program   &  $51,000$  & $300$ &  $2$ & $200$ & $1024$ & Adam &$512$ & $96.47\%$  & $96.42\%$     \tabularnewline 
					\bottomrule
			\end{tabular}}
		}
	\end{center}
	\label{tab:classifier parameters}
\end{table*}

We modified the New England 39-bus system case to include $4$ wind generators connected to buses $5,\,6,\,14$ and $17$ and used the scenario-based stochastic method presented in \cite{ETH-AC} to formulate the sampled optimisation problem. The number of uncertain parameters in the problem is 4---corresponding to uncertain active power generated by renewable generators. The penetration level of renewable generators is $30\%$, meaning that renewable generators can provide up to $30\%$ of the total demand. The uncertainty distribution is chosen by the Pearson system with a standard deviation equal to $0.2\times$ (predicted generation power) and kurtosis of $3.5$ leading to a leptokurtic distribution with a heavier tail than that of a Gaussian. 
The configuration and parameters for training the neural network classifier are tabulated in Table \ref{tab:classifier parameters}. 
In Table \ref{tab:OPF results}, we report the time that it takes to directly solve the sampled problem using Mosek \cite{andersen2000mosek}, Algorithm \ref{alg:sequential algorithm}, and Algorithm \ref{alg:learning sequential algorithm} for different number of scenario samples. For small number of scenario samples, e.g. $N =100$, Algorithm \ref{alg:sequential algorithm} is slower than directly solving the sampled problem using Mosek, however, Algorithm\ref{alg:learning sequential algorithm} is still faster than Mosek. For large number of samples, Algorithms \ref{alg:sequential algorithm} and \ref{alg:learning sequential algorithm} both outperform Mosek. For instance, when the number of scenario samples is  $10^4$, Algorithms \ref{alg:sequential algorithm} and \ref{alg:learning sequential algorithm} are respectively $21$, and $32$ times faster than Mosek. 

\begin{table}[t]
	\caption{CPU time taken to directly solve the sampled optimal power flow problem using Mosek, and the CPU time it takes for Algorithms \ref{alg:sequential algorithm} and \ref{alg:learning sequential algorithm} to solve the problem for different values of the scenario samples.}
	\begin{center}
		\begin{tabular}{c||c|c|c}
			\toprule
			\# Scenario & CPU Time    & CPU Time  & CPU Time  \tabularnewline
			
			Samples & Mosek  &  Algorithm \ref{alg:sequential algorithm} & Algorithm \ref{alg:learning sequential algorithm} \tabularnewline
			\midrule
			\midrule
			$10^2$ & $121.4$  &  $127.9$  & $56.3$      \tabularnewline \midrule
			$10^3$ & $1598.8$  &  $363.9$  & $165.7$      \tabularnewline \midrule
			$5\times 10^3$ & $3.06\times 10^{4}$  &  $2594.4$  & $1286.7$      \tabularnewline \midrule
			$10^4$ & $1.18\times 10^5$  &  $5657.9$  & $3686.5$      \tabularnewline

			\bottomrule
		\end{tabular}	
	\end{center}
	\label{tab:OPF results}
\end{table}

\subsection{Robust Unit Commitment}
Unit commitment is a mathematical optimisation problem solved in power grids to determine the commitment of each generator. It considers a time horizon, and given a predicted demand over the considered horizon and generators' minimum and maximum power ratings, its solution defines which generators should be online and which ones should be offline so that the total generation cost is minimised. There are many models developed for the unit commitment problem in the literature, see \cite{haaberg2019fundamentals} and references therein for a full review of the topic. In this subsection, we consider a simplistic version of this problem. Sets, indices, and variables used in defining the model are introduced first. 

\begin{tabular}{p{0.2\columnwidth} p{0.75 \columnwidth}}
	$\mathcal{G}$& set of generators in the grid with cardinality $n_g,\, |\mathcal{G}|=n_g$\\
	$T\in\mathbb{N}$ & time horizon over which the problem is solved \\
	$t\in\{1,\ldots,T\}$ & time periods\\
	$P_{i,t}\in\mathbb{R}$ & active power generated by generator $i$ at time period $t$\\
	$U_{i,t}\in\{0,1\}$ & on-off status of generator $i$ at time $t$\\
	$D_t\in \mathbb{R}$ & demand at time $t$\\
	$P^{\min}_{i,t}\in\mathbb{R}$ & minimum active power generator $i$ should provide at time $t$ \\ 
	$P^{\max}_{i,t}\in\mathbb{R}$ & maximum active power generator $i$ can provide at time $t$ \\
	$\tau^\text{OFF}_i$ &  unit $i$ must be off-line for $\tau^\text{OFF}_i$ before it can be on-line\\  
	$\tau^\text{ON}_i$ &  unit $i$ must be on-line for $\tau^\text{ON}_i$ before it can be off-line\\  
	$\Delta P^{\max}_{i,t}$ &   maximum allowed difference between power generated by generator $i$ at time $t$ and $t-1$
\end{tabular}	
\textbf{Objective}\\
The objective is to minimise the total operating cost of all generators across the grid
\[
f(P) = \sum_{i}^{n_g}\sum_{t=1}^{T} Q_{ii}P^2_{i,t}+C_iP_{i,t},
\] 
where $Q\in\mathbb{R}^{n_g\times n_g}$ and $C\in\mathbb{R}^{n_g}$ are, respectively, diagonal matrix and vector defining the running cost of generators. 

\textbf{Constraints}\\
The amount of power each generator can provide is constrained by the following constraint
\begin{equation*}
U_{i,t}P^{\min}_{i,t}\leq P_{i,t}\leq U_{i,t}P^{\max}_{i,t}.
\end{equation*}
The total power generated by active generators should meet demand at all time
\begin{equation*}
\sum_{i=1}^{n_g}P_{i,t}\geq D_t, \,\, \forall t=1,\ldots,T.
\end{equation*}
The minimum up-time and downtime of each generator are defined using the following constraint
\begin{equation*}
U_{i,\tau}\geq U_{i,t}-U_{i,t-1},\, \tau = t,t+1,\ldots,\min(T,t+\tau^\text{ON}_i-1),\,\forall i\in\mathcal{G},\,t = 2,\ldots,T.
\end{equation*} 
\begin{equation*}
U_{i,\tau}\leq 1- U_{i,t-1}-U_{i,t},\, \tau = t,t+1,\ldots,\min(T,t+\tau^\text{OFF}_i-1),\,\forall i\in\mathcal{G},\,t = 2,\ldots,T.
\end{equation*} 
The above two constraints require generator $i$ to remain online (resp. offline) for $\tau^\text{ON}_i$ (resp. $\tau^\text{OFF}_i$) time periods before they can go offline (resp. online). 
The following ramp constraint limits the rate of change of active power generated by each generator at each sampling time
\begin{equation*}
P_{i,t}-P_{i,t-1}\leq\Delta P_{i,t}^{\max}, \, \forall t = 1,\ldots,T,\, \forall i\in\mathcal{G}.
\end{equation*}
The demand at time $t$, denoted as $D_t$, is not fully known. To capture the uncertainty associated with demand, we assume that $D_t$ is constructed by a nominal predicted demand $D^0_t$ and an uncertain demand $D^q_t$:
\[
D_t = D^0_t + D^q_t, \,\,t = 1,\ldots,T.
\] 
\begin{table}[t]
\caption{CPU time it takes to directly solve the sampled unit commitment problem using Gurobi, and the CPU time it takes for Algorithms \ref{alg:sequential algorithm} and \ref{alg:learning sequential algorithm} to solve the problem for different values of the scenario samples.}
\begin{center}
	\begin{tabular}{c||c|c|c}
		\toprule
		\# Scenario & CPU Time    & CPU Time  & CPU Time  \tabularnewline
		
		Samples & Gurobi  &  Algorithm \ref{alg:sequential algorithm} & Algorithm \ref{alg:learning sequential algorithm} \tabularnewline
		\midrule
		\midrule
		$10^2$ & $17.2$  &  $27.1$  & $11.6$      \tabularnewline \midrule
		$10^3$ & $288.2$  &  $36.1$  & $14.5$      \tabularnewline \midrule
		$5\times 10^3$ & $3904.1$  &  $68.7$  & $42.2$      \tabularnewline \midrule
		$10^4$ & $2.3\times 10^4$  &  $106.9$  & $72.9$      \tabularnewline 			
		\bottomrule
	\end{tabular}	
\end{center}
\label{tab:unit commitment results}
\end{table}
For computational purposes, we select $n_g=4,\,T=12,\,P^{\max}_{i,t} = \texttt{randint(1,115)},\,P^{\min}_{i,t}=\lceil P^{\max}_{i,t}/2\rceil, \, D_t^0=150\sin(2\pi/24\,t),t=1,\ldots,T,\,Q=\text{diag}(\texttt{randint(1,50)}),\,C = \texttt{randint(1,50)},\, \tau^{\text{ON}}_i=\texttt{randint(1,T)},\,\tau^{\text{OFF}}_i=\texttt{randint(1,T)},\forall i \in \mathcal{G}$. The uncertain component of the demand is bounded in $[-1,1]$, i.e. $D_t^q\in[-1,1], \forall t = 1,\ldots,T$. The configuration and parameters used for training the neural network classifier used in Algorithm \ref{alg:learning sequential algorithm} is shown in Table  \ref{tab:classifier parameters}. We used Gurobi version $9.1.2$ to solve the sampled optimisation problem for different number of scenario samples and compared its performance in terms of the time it takes to solve the problem with Algorithms \ref{alg:sequential algorithm}, and \ref{alg:learning sequential algorithm}. The result of this simulation is shown in Table \ref{tab:unit commitment results}. Similar to Table \ref{tab:OPF results}, for a small number of samples directly solving the sampled optimisation problem using Gurobi  results in a shorter solution time than Algorithm \ref{alg:sequential algorithm}. However, for a large number of scenarios, both Algorithms \ref{alg:sequential algorithm} and \ref{alg:learning sequential algorithm} notably outperform Gurobi. For $N=10,000$ (last row of Table \ref{tab:unit commitment results}), Algorithms \ref{alg:sequential algorithm}  and \ref{alg:learning sequential algorithm} are respectively $217$ and $319$ times faster than the direct solution obtained using Gurobi. This shows the significant computational improvement one can achieve by using Algorithms \ref{alg:sequential algorithm} and \ref{alg:learning sequential algorithm}.

\subsection{Robust Mixed-integer Linear Programs}\label{sec:lp}
Classical robust Mixed-Integer Linear Programming (MILP) problems admit the following formulation
\begin{align}\label{eq:milp}
&\text{minimize} \qquad c^T\decision\\ \nonumber
&\text{subject to: } A\decision\leq b+b_q,  \\ \nonumber
& \qquad\qquad\decision\in\bR^{25}\times \bZ^{5},
\end{align} 
where the objective definition is defined by the vector $c\in\bR^{30}$, while $A\in\bR^{500\times 30},\, b\in\bR^{500}$ constitute the (fixed) matrix and vector used to define the set of nominal constraints of Problem \eqref{eq:milp}, and $b_q\in\bR^{500}$ is a so-called interval vector---i.e., a vector whose entries vary in given intervals---characterising the uncertainty in the optimisation problem \eqref{eq:milp}.  The vectors $b,c$ and nominal matrix $A$ are generated such that problem \eqref{eq:milp} is feasible. 
To this end, we follow the methodology presented in
\citep{dunham_experimental_1977}.
%
%
The distribution of uncertain vector $b_q$ is uniform and its entries are bounded in $b\times[-0.01, 0.01]$.
%
%
The sampled version of problem \eqref{eq:milp} is constructed by extracting random samples $\{b_q^{(i)}\}_{i=1}^N$ from the set of uncertainty  
\begin{align}\label{eq:sampled milp}
&\text{minimize} \qquad c^T\decision\\ \nonumber
&\text{subject to: } A\decision\leq b+b_q^{(i)},  \,\, i=1,\ldots,N\\ \nonumber
& \qquad\qquad\decision\in\bR^{25}\times \bZ^{5}.
\end{align} 
The hyper-parameters used in training the deep neural network classifier are listed in Table \ref{tab:classifier parameters}.
In Table \ref{tab: simulation results MILP}, we vary the number of scenario samples $N$ and solve problem \eqref{eq:sampled milp} using Algorithms \ref{alg:sequential algorithm} and, \ref{alg:learning sequential algorithm} and compare their performance against directly solving Problem \eqref{eq:sampled milp} using Gurobi \cite{gurobi}.  
%
For a small number of scenario samples, Algorithm \ref{alg:sequential algorithm} is slower than directly solving the sampled optimisation problem \eqref{eq:sampled milp} using Gurobi. However, for all the scenario samples, Algorithm \ref{alg:learning sequential algorithm} outperforms Gurobi and Algorithm \ref{alg:sequential algorithm}. It is worth noting that for some entries in Table \ref{tab: simulation results MILP}, the Gurobi solver has been found to require more than $64$ GB of RAM to complete the solution process, thus preventing it from completing this task on the cluster. This highlights yet another major advantage of Algorithms \ref{alg:sequential algorithm} and \ref{alg:learning sequential algorithm} in the fact that they achieve significant memory savings compared to classical algorithms meant to solve sampled optimisation problems associated with a large number of constraints. 

\begin{table}[t]
\caption{CPU time it takes to directly solve the sampled MILP problem using Gurobi, and the CPU time it takes for Algorithms \ref{alg:sequential algorithm} and \ref{alg:learning sequential algorithm} to solve the problem for different values of the scenario samples. NA refers to the case that Gurobi requires more resources to solve the problem.}
\begin{center}
	\begin{tabular}{c||c|c|c}
		\toprule
		\# Scenario & CPU Time    & CPU Time  & CPU Time  \tabularnewline
		
		Samples & Gurobi  &  Algorithm \ref{alg:sequential algorithm} & Algorithm \ref{alg:learning sequential algorithm} \tabularnewline
		\midrule
		\midrule
		$10^4$ & $89.3$  &  $232$  & $78.8$      \tabularnewline \midrule			
		$5\times 10^4$ & $466.5$  &  $292.2$  & $109.9$      \tabularnewline \midrule
		 $10^5$ & $1220$  &  $295.2$  & $161.1$      \tabularnewline \midrule
		 $5\times 10^5$ & NA  &  $615.5$  & $452.1$      \tabularnewline 
		\bottomrule
	\end{tabular}	
\end{center}
\label{tab: simulation results MILP}
\end{table}

\section{Conclusion}
In this paper, we presented two algorithms for solving sampled optimisation problems. Both algorithms exhibit a significant saving in time and memory required for solving this class of optimisation problems. 
Both algorithms involve two main steps: verification and optimisation, which are performed sequentially to converge toward the optimal solution. At each step of these algorithms, we need to compute the basis--- a minimal set of constraints defining the current solution. Algorithm \ref{alg:learning sequential algorithm} features a neural network multiclass classifier to reduce the complexity associated with finding the basis at each iteration of the algorithm. 
The convergence properties of both algorithms are analysed, and extensive numerical simulations are performed to compare their performance---in solving various non-trivial sampled optimisation problems---against widely used commercial solvers.

The two proposed algorithms significantly reduce the computational time of solving problems for which the number of constraints is much larger than the number of decision variables. If the number of decision variables is large, the combinatorial dimension of the problem might grow---see Theorem \ref{thm:combinatorial dim} for the exact upper bound---leading to a possible increase in the complexity of Algorithms \ref{alg:sequential algorithm} and \ref{alg:learning sequential algorithm}. A  possible future direction is to combine the sequential nature of the two proposed algorithms with column generation methods \cite{desaulniers2006column,ford1958suggested} to reduce the computational complexity for the case that the number of decision variables is large.



\end{document}